\documentclass[12pt,a4paper]{article}
\usepackage{amsfonts,amssymb,amsmath,amsthm,indentfirst}

\numberwithin{equation}{section}

\begin{document}

\date{}
\title{Hardy-Littlewood-Sobolev Type Inequality and Stein-Wiess Type Inequality on Carnot Groups}

\author{Tingxi Hu;
Pengcheng Niu*\\
\small\it Department of Applied Mathematics; Key Laboratory of Space
Applied Physics and Chemistry,\\
\small\it Ministry of Education, Northwestern Polytechnical
University, Xi'an, Shaanxi,\\
\small\it710129, P. R. China
 }
\date{}
\maketitle


A Stein-Weiss type inequality on Carnot groups is established by proving the
boundedness of an integral operator and the Hardy-Littlewood-Sobolev type
inequality on Carnot groups is also derived.


\section{Introduction}
The Hardy-Littlewood-Sobolev inequality on $R^N$ (i.e., HLS
inequality, see [1], [2] and [3]) is of the form
\begin{equation}
\left| {\iint_{R^N  \times R^N } {\frac{{\overline {f\left( x
\right)} g\left( y \right)}} {{\left| {x - y} \right|^\lambda
}}dxdy}} \right| \leqslant C_{r,\lambda ,N} \left\| f \right\|_r
\left\| g \right\|_s,
\end{equation}
where $1 < r,s < \infty ,0 < \lambda  < N$ and $\frac{1} {r} +
\frac{1} {s} + \frac{\lambda } {N} = 2$, $C_{r,\lambda ,N}$ is a
positive constant. Stein and Wiess in [4] proved a weighted HLS
inequality on $R^N$, known as the Stein-Wiess inequality, which
reads
\begin{equation}
\left| {\iint_{R^N  \times R^N } {\frac{{\overline {f\left( x
\right)} g\left( y \right)}} {{\left| x \right|^\alpha  \left| {x -
y} \right|^\lambda  \left| y \right|^\beta  }}dxdy}} \right|
\leqslant C_{\alpha ,\beta ,r,\lambda ,N} \left\| f \right\|_r
\left\| g \right\|_s,
\end{equation}
where $1 < r,s < \infty ,0 < \lambda  < N, \frac{1} {r} + \frac{1}
{s} + \frac{\lambda } {N} = 2$ and $\alpha+\beta\geq0$, such that
$\lambda+\alpha+\beta\leq N$, $\alpha  < N/r'$, $\beta  < N/s'$,
here $r'$ and $s'$ denote respectively the conjugate indice of $r$
and $s$, $C_{\alpha,\beta,r,\lambda,N}$ is a positive constant
independent of the functions $f,g$.

The HLS type inequality on the Heisenberg group was proved by
Folland and Stein see [5, Proposition 8.7 and Lemma 15.3].
Recently, Han, Lu and Zhu in [6] extended (1.1) to the Heisenberg
group and derived Stein-Wiess type inequalities with weights.

In this paper we pay attention to Stein-Wiess type inequality on
Carnot groups. As a consequence, the HLS type inequality on Carnot
groups is derived. The proofs of these results are motivated by ones in [6].

This paper is organized as follows. In Section 2 we recall some
basic facts about Carnot groups and state precisely our assumptions
and main results; in Section 3 and Section 4 we prove Theorems 2.1
and 2.2, respectively. Section 5 is devoted to the proof of Lemma
4.4 used in Section 4.

\section{Carnot groups and main results}
We begin by describing Carnot groups. For more information, we
refer to [7] and references therein. A Carnot group $G$ of step $r$ is
a simply connected nilpotent Lie group such that its Lie algebra
$\mathfrak{g}$ admits a stratification
$$\mathfrak{g}= {V_1} \oplus {V_2} \oplus  \ldots  \oplus {V_r} =  \oplus _{l = 1}^rV,$$
with $[{V_1},{V_l}] = {V_{l + 1}}\;\left( {l = 1,2, \ldots ,r - 1}
\right)$ and $[{V_1},{V_l}] = \left\{ 0 \right\}$.

Denoting ${m_l} = \dim {V_l}$, we fix on $G$ a system of coordinates
$u = \left( {{z_1},{z_2}, \ldots ,{z_r}} \right)$, ${z_l} \in {R^{{m_l}}}$.
Every Carnot group $G$ is naturally equipped with a family of non-isotropic
dilations defined by ${\delta _r}(r > 0)$:%
$${\delta _r}\left( u \right) = \left( {r{z_1},{r^2}{z_2}, \ldots ,{r^r}{z_r}}
 \right)\; ,\;for\;u \in G;$$
the homogeneous dimension of $G$ is given by $Q = \sum\limits_{l =
1}^r {l{m_l}} $. We denote by $du$ a fixed Haar measure on $G$. One
easily has $\left( {d \circ {\delta _r}} \right)\left( u \right) =
{r^Q}du$.

The group law given by Baker-Campbell-Hausdorff formula is%
$$uv = u + v + \sum\limits_{1 \leqslant l,k \leqslant r} {{Z_{l,k}}\left( {u,v} \right)}, \;for \;u,v \in G,$$
where each $Z_{l,k} \left( {u,v} \right)$ is a fixed linear
combination of iterated commutators containing $l$ times $u$ and $k$
times $v$. The homogenous norm of $u$ on $G$ is defined by%
$$\left| u \right| = \left( {\sum\limits_{j = 1}^r {\left| {z_j }
\right|^{\frac{{2r!}} {j}} } } \right)^{\frac{1} {{2r!}}},$$ where
$\left| {z_j } \right|$ denotes the Euclidean distance $z_j  \in
R^{m_j }$ to the origin. Such homogenous norm on $G$ can be used to
define
a pseudo-distance on $G$:%
$$d\left( {u,v} \right) = \left| {u^{ - 1} v} \right|.$$Denote the pseudo-ball of
radius $r$ about $u$ by $B\left( {u,r} \right) = \left\{ {v \in
G|d\left( {u,v} \right) < r} \right\},$ and the unit pseudo-ball
about origin by $\left\{ {\left| u \right| < 1} \right\}.$

We introduce $L^{p}$ norm with the weight $w$ on $G$ with:%
$$\left\| f \right\|_{L^p \left( {G,w} \right)}  = \left( {\int_G
{\left| {f\left( u \right)} \right|^p wdu} } \right)^{\frac{1}
{p}},$$ and denote the space of all measure functions with finite
weighted $L^{p}$ norm by $L^p \left( {G,w} \right).$

Our main results in this paper are

\textbf{Theorem 2.1 (HLS Inequality with Weight $\left| u \right|$ ,
i.e., Stein-Wiess type Inequality).} Let $1 < r,s < \infty ,0 < \lambda  < Q$ and $\alpha  + \beta
\geqslant 0$ such that $\lambda  + \alpha  + \beta  \leqslant
Q,\alpha  < Q/r',\beta  < Q/s'$ and $\frac{1} {r} + \frac{1} {s} +
\frac{{\lambda  + \alpha  + \beta }} {Q} = 2$, then there exists a
positive constant $C_{\alpha ,\beta ,r,\lambda ,G}$ independent of
the functions $f,g$ such that%
\begin{equation}
\left| {\iint_{G \times G} {\frac{{\overline {f\left( u \right)}
g\left( v \right)dudv}} {{\left| u \right|^\alpha  \left| {u^{ - 1}
v} \right|^\lambda  \left| v \right|^\beta  }}}} \right| \leqslant
C_{\alpha ,\beta ,r,\lambda ,G} \left\| f \right\|_r \left\| g
\right\|_s ,
\end{equation}
where $u = (z_1 ,z_2 , \ldots ,z_r ),v = (z'_1 ,z'_2 , \ldots
,z'_r ) \in G.$

\textbf{Theorem 2.2 (HLS Inequality with Weights $\left| {z_l }
\right|\left( {1 \leqslant l \leqslant r} \right)$).} Let $1 < r,s <
\infty ,0 < \lambda  < Q$ and $0 \leqslant \alpha  + \beta  <
\frac{{m_l }} {{Q - lm_l }}\lambda$ such that $\lambda  + l\alpha  +
l\beta  \leqslant Q,\alpha  < m_l /r',\beta  < m_l /s'$ and
$\frac{1} {r} + \frac{1} {s} +\frac{{\lambda  + l\alpha  + l\beta }}
{Q} = 2$, then there exists a positive constant $C_{\alpha ,\beta
,r,\lambda ,G,l}$ independent of the functions $f,g$ such that%
\begin{equation}
\left| {\iint_{G \times G} {\frac{{\overline {f\left( u \right)}
g\left( v \right)dudv}} {{\left| {z_l } \right|^\alpha  \left| {u^{
- 1} v} \right|^\lambda  \left| {z'_l } \right|^\beta  }}}} \right|
\leqslant C_{\alpha ,\beta ,r,\lambda ,G,l} \left\| f \right\|_r
\left\| g \right\|_s,
\end{equation}
where $u = (z_1 ,z_2 , \ldots ,z_r ),v = (z'_1 ,z'_2 , \ldots
,z'_r ) \in G.$

From Theorem 2.1 or Theorem 2.2 we immediately obtain

\textbf{Corollary 2.3 (HLS type Inequality on Carnot groups).}
Let $1 < r,s < \infty ,0 < \lambda  < Q$ such that $\frac{1} {r} +
\frac{1} {s} + \frac{\lambda } {Q} = 2,$ then there exists a
positive constant $C_{r,\lambda ,G}$ independent of the functions
$f,g$ such that%
\begin{equation*}
\left| {\iint_{G \times G} {\frac{{\overline {f\left( u \right)}
g\left( v \right)}} {{\left| {u^{ - 1} v} \right|^\lambda  }}dudv}}
\right| \leqslant C_{r,\lambda ,G} \left\| f \right\|_r \left\| g
\right\|_s,
\end{equation*}
where $u = (z_1 ,z_2 , \ldots ,z_r ),v = (z'_1 ,z'_2 , \ldots
,z'_r ) \in G.$

\section{Proof of Theorem 2.1}
We first state two versions which are actually equivalent to
Theorem 2.1 (see Proposition 3.3 bellow).

\textbf{Theorem 3.1.} Let $1 < p \leqslant q < \infty ,0 < \lambda <
Q$ and $\alpha  + \beta  \geqslant 0$ such that $\alpha  < Q/q,
\beta  < Q/p'$ and $\frac{1} {q} = \frac{1} {p}+\frac{{\lambda  + \alpha  + \beta
}} {Q} - 1,$ then%
\begin{equation}
\left\| {Tf} \right\|_{L^q (G,\left| u \right|^{ - \alpha q} )}
\leqslant C\left\| f \right\|_{L^p (G,\left| u \right|^{\beta p} )},
\end{equation}
where $T$ is an integral operator with the form%
\begin{equation}
Tf\left( u \right) := T_\lambda  f\left( u \right) = \int_G
{\frac{{f\left( v \right)dv}} {{\left| {u^{ - 1} v} \right|^\lambda
}}},
\end{equation}
$C = C_{\alpha ,\beta ,p,\lambda ,G,l}$ is a positive constant
independent of $f$.

We notice that (3.1) means%
\begin{equation}
\left( {\int_G {\left| {Tf\left( u \right)} \right|^q \left| u
\right|^{ - \alpha q} du} } \right)^{\frac{1} {q}}  \leqslant
C\left( {\int_G {\left| {f\left( u \right)} \right|^p \left| u
\right|^{\beta p} du} } \right)^{\frac{1} {p}},
\end{equation}
and easily get a equivalent following form of Theorem 3.1 by choosing $
f\left( v \right) = \frac{{g\left( v \right)}} {{\left| v
\right|^\beta  }}$:

\textbf{Theorem 3.2.} Let $1 < p \leqslant q < \infty ,0 < \lambda <
Q$ and $\alpha  + \beta  \geqslant 0$ such that $\alpha  < Q/q,
\beta  < Q/p'$ and $\frac{1} {q} =\frac{1} {p}+ \frac{{\lambda  + \alpha  + \beta
}} {Q} - 1$, then%
\begin{equation*}
\left\| {Sg} \right\|_q  \leqslant C\left\| g \right\|_p,
\end{equation*}
where $S$ is an integral operator with%
\begin{equation*}
Sg\left( u \right) := S_{\lambda ,\alpha ,\beta } g\left( u \right) =
\int_G {\frac{{g\left( v \right)dv}} {{\left| u \right|^\alpha
\left| {u^{ - 1} v} \right|^\lambda  \left| v \right|^\beta  }}}
\end{equation*}
and $C = C_{\alpha ,\beta ,p,\lambda ,G,l}$ is a positive constant
independent of $g$.

If we denote $s = p$ and $r = q'$, then we see with some simple
computations that parameters in Theorems 2.1, 3.1 and 3.2 are same.

\textbf{Proposition 3.3.} Theorems 2.1 and 3.2 are equivalent.

\textbf{Proof.} If Theorem 2.1 is true, then, it follows by the dual
argument that%
\begin{eqnarray*}
\left\| {Sg} \right\|_q &\leq&\mathop {\sup }\limits_{\left\| f
\right\|_{q'}  = 1} \left| {\int_G {\overline {f\left( u \right)}
Sg\left( u \right)du} } \right| \\
&=& \mathop {\sup }\limits_{\left\| f \right\|_{q'}  = 1} \left|
{\iint_{G \times G} {\frac{{\overline {f\left( u \right)} g\left( v
\right)dudv}} {{\left| u \right|^\alpha  \left| {u^{ - 1} v}
\right|^\lambda  \left| v \right|^\beta  }}}} \right| \\
&\leq& C\mathop {\sup }\limits_{\left\| f \right\|_{q'}  = 1}
\left\| f \right\|_{q'} \left\| g \right\|_p\\
&=&C\left\| g \right\|_p .
\end{eqnarray*}
On the other hand, we compute by using the inequality in Theorem 3.2 to have%
\begin{eqnarray*}
\left| {\iint_{G \times G} {\frac{{\overline {f\left( u \right)}
g\left( v \right)dudv}} {{\left| u \right|^\alpha  \left| {u^{ - 1}
v} \right|^\lambda  \left| v \right|^\beta  }}}} \right| &=&
\left|{\int_G {\overline {f\left( u \right)} \left( {\int_G
{\frac{{g\left( v \right)dv}} {{\left| u \right|^\alpha  \left| {u^{
- 1} v} \right|^\lambda  \left| v \right|^\beta  }}} } \right)du} }
\right| \\
&=& \left| {\int_G {\overline {f\left( u \right)} Sg\left( u
\right)du}} \right|\\
&\leq&\left\| f \right\|_r \left\| {Sg} \right\|_{r'} \\
&\leq&C\left\| f \right\|_r \left\| g \right\|_s .
\end{eqnarray*}
The proof in ended.$\Box$

From Proposition 3.3 and the equivalence between Theorems 3.1 and
3.2, we immediately get Theorem 2.1 if one proves Theorem 3.1. To do
this, we need some preparations.

Denoting $K\left( {u,v} \right) = \left| {u^{ - 1} v} \right|^{ -
\lambda }  = d\left( {u,v} \right)^{ - \lambda }$, it see that (3.2) is changed
into%
\begin{equation*}
Tf\left( u \right) = \int_G {K\left( {u,v} \right)f\left( v
\right)dv}.\eqno(3.2')
\end{equation*}
For the application later, let us note that $d\left( {u,v} \right)$ satisfies the triangle
inequality%
\begin{equation}
d\left( {u_1 ,u_2 } \right) \leqslant K_G \left[ {d\left( {u_1 ,u_3
} \right) + d\left( {u_3 ,u_2 } \right)} \right].
\end{equation}

A positive measure $wdu$ is said a doubling measure on $G$ if there
exists some positive constant $C$ such that%
\begin{equation*}
\int_B {wdu}  \leqslant C\int_{B'} {wdu}
\end{equation*}
for all pair of concentric ball $B$ and $B'$ satisfying $r\left( B
\right) = 2r\left( {B'} \right)$. Clearly, since%
\begin{equation*}
\int_{\left| u \right| < 2r} {\left| u \right|^{ - \alpha p\tau }
du}  = \int_{\left| {\delta _2  \circ u} \right| < 2r} {\left|
{\delta _2  \circ u} \right|^{ - \alpha p\tau } \left( {d \circ
\delta _2 } \right)\left( u \right)}  = 2^{ - \alpha p\tau  + Q}
\int_{\left| u \right| < r} {\left| u \right|^{ - \alpha p\tau }
du},
\end{equation*}
we find that $\left| u \right|^{ - \alpha q\tau } du$ is a doubling
measure on $G$. Similarly, $\left| u \right|^{\beta p\tau } du,
\left| {z_l } \right|^{ - \alpha q\tau } du, \left| {z_l }
\right|^{\beta p\tau } du$ are also doubling measures on $G$

Sawyer and Wheeden in [8, p.821] demonstrated boundedness for a
class of integral operators on the homogeneous space. Noting that
the Carnot group is a special homogeneous space, we can apply Sawyer
and Wheeden's result to the Carnot group. For convenience, we
restate the result in [8] on the Carnot group.

\textbf{Lemma 3.4.} Let $w_1 \left( u \right)$ and $w_2 \left( u
\right)$ are nonnegative functions on $G$. Then the operator $T$ in
$(3.2')$ is bounded from $L^p \left( {G,w_2 } \right)$ to $L^q
\left( {G,w_1 } \right)$if both of the following two statements are
true:

(1). There exists $\varepsilon  > 0$ such that for any pair of balls
$B$ and $B'$ with radius $r$ and $r'$ satisfying $B' \subseteq 4B$,
it holds%
\begin{equation}
\left( {\frac{{r'}} {r}} \right)^{Q - \varepsilon} \left(
{\frac{{\varphi \left( {B'} \right)}} {{\varphi \left( B \right)}}}
\right) \leqslant C_\varepsilon,
\end{equation}
where $\varphi \left( B \right) = \sup \left\{ {K\left( {u,v}
\right)|u,v \in B,d\left( {u,v} \right) \geqslant C\left( {K_G }
\right)r} \right\}$ for a ball $B\subseteq G$ with radius $r$, $
C\left( {K_G } \right) = K_G^{ - 4} /9$, $C_{\varepsilon}$ is a
positive constant depending only on $\varepsilon$.

(2). There exists $\tau  > 1$ such that $w_1^\tau  du$ and $
w_2^{\left( {1 - p'} \right)\tau } du$ are doubling measures, and%
\begin{equation}
\varphi \left( B \right)\left| B \right|^{\frac{1} {{p'}} + \frac{1}
{q}} \left( {\frac{1} {{\left| B \right|}}\int_B {w_1^\tau  du} }
\right)^{\frac{1} {{q\tau }}} \left( {\frac{1} {{\left| B
\right|}}\int_B {w_2^{\left( {1 - p'} \right)\tau } du} }
\right)^{\frac{1} {{p'\tau }}}  \leqslant C_\tau,
\end{equation}
for all ball $B \subseteq G$, where $C_\tau$ is a positive constant
depending only on $\tau$.

\textbf{Proof of Theorem 3.1.} We will apply the above Lemma 3.4 to
verify Theorem 3.1, that is, we will show (3.5) and (3.6)
in Lemma 3.4 are satisfied for the operator defined in $(3.2')$ by
choosing $ w_1 \left( u \right) = \left| u \right|^{ - \alpha q}$
and $w_2 \left( u \right) = \left| u \right|^{\beta p}$.

In fact, since $K\left( {u,v} \right) = \left| {u^{ - 1} v}
\right|^{ - \lambda }  = d\left( {u,v} \right)^{ - \lambda }$, it
obtains $ \varphi \left( B \right) = 9^\lambda  K_G^{4\lambda } r^{
- \lambda } ,\varphi \left( {B'} \right) = 9^\lambda K_G^{4\lambda }
\left( {r'} \right)^{ - \lambda }$ and%
\begin{equation*}
\left( {\frac{{r'}} {r}} \right)^Q \left( {\frac{{\varphi \left(
{B'} \right)}} {{\varphi \left( B \right)}}} \right) = \left(
{\frac{{r'}} {r}} \right)^{Q - \lambda }.
\end{equation*}
Due to $0 < \lambda  < Q$, we let $\varepsilon  > 0$ small such that
$ Q - \lambda  - \varepsilon  > 0$, and see that for $B' \subseteq
4B$,%
\begin{equation*}
\left( {\frac{{r'}} {r}} \right)^{Q - \varepsilon } \left(
{\frac{{\varphi \left( {B'} \right)}} {{\varphi \left( B \right)}}}
\right) = \left( {\frac{{r'}} {r}} \right)^{Q - \lambda  -
\varepsilon }  \leqslant 4^{Q - \lambda  - \varepsilon }  =
C_\varepsilon   < \infty,
\end{equation*}
which proves (3.5).

Furthermore, denoting%
\begin{eqnarray*}
M_1 &=&\varphi \left( B \right)\left| B \right|^{\frac{1} {{p'}} +
\frac{1} {q}},\\
M_2 &=&\left( {\frac{1} {{\left| B \right|}}\int_B {\left| u
\right|^{ - \alpha q\tau } du} } \right)^{\frac{1} {{q\tau }}},\\
M_3 &=&\left( {\frac{1} {{\left| B \right|}}\int_B {\left| u
\right|^{ - \beta p'\tau } du} } \right)^{\frac{1} {{p'\tau }}},
\end{eqnarray*}
we estimate these $M_i \;(i = 1,2,3)$. Let $\bar \lambda {\text{ =
}}Q\left( {\frac{1} {q} + \frac{1} {{p'}}} \right)$, then $ \frac{1}
{q} = \frac{1} {p} + \frac{{\bar \lambda }} {Q} - 1$ and $\bar
\lambda  - \lambda  - \alpha  - \beta  = 0$. A direct computation
implies%
\begin{equation*}
M_1  = \varphi \left( B \right)\left| B \right|^{\frac{1} {q} +
\frac{1} {{p'}}}  = Cr^{ - \lambda } r^{Q\left( {\frac{1} {q} +
\frac{1} {{p'}}} \right)}  = Cr^{\bar \lambda  - \lambda };
\end{equation*}
if $\alpha q\tau  < Q$,%
\begin{equation*}
M_2  = \left( {\frac{1} {{\left| B \right|}}\int_B {\left| u
\right|^{ - \alpha q\tau } du} } \right)^{\frac{1} {{q\tau }}}
\leqslant \left( {C_\tau  \frac{{r^{Q - \alpha q\tau } }} {{r^Q }}}
\right)^{\frac{1} {{q\tau }}}  = C_\tau  r^{ - \alpha };
\end{equation*}
and if $\beta p'\tau  < Q$,%
\begin{equation*}
M_3  = \left( {\frac{1} {{\left| B \right|}}\int_B {\left| u
\right|^{ - \beta p'\tau } du} } \right)^{\frac{1} {{p'\tau }}}
\leqslant \left( {C_\tau  \frac{{r^{Q - \beta p'\tau } }} {{r^Q }}}
\right)^{\frac{1} {{p'\tau }}}  = C_\tau  r^{ - \beta }.
\end{equation*}
Since $\alpha  < Q/q,\beta  < Q/p',$ it follows $\min \left\{
{\frac{Q} {{\alpha q}},\frac{Q} {{\beta p'}}} \right\} > 1$ which
ensures existence of $\tau$ such that $1 < \tau  < \min \left\{
{\frac{Q} {{\alpha q}},\frac{Q} {{\beta p'}}} \right\}$. Thus $
\alpha q\tau  < Q,\beta p'\tau  < Q,$ and%
\begin{equation*}
M_1  \cdot M_2  \cdot M_3  \leqslant C_\tau  r^{\bar \lambda  -
\lambda  - \alpha  - \beta }  = C_\tau   < \infty,
\end{equation*}
which shows (3.6). Then by Lemma 3.4, the proof of Theorem 3.1 is
completed.$\Box$

\section{Proof of Theorem 2.2}
We need the following conclusion to prove Theorem 2.2.

\textbf{Theorem 4.1.} Let $1 < p \leqslant q < \infty ,0 < \lambda <
Q,1 \leqslant l \leqslant r$ and $ 0 \leqslant \alpha  + \beta  <
\frac{{m_l }} {{Q - lm_l }}\lambda$ such that $\alpha  < m_l /q,
\beta  < m_l /p',\frac{1} {q} = \frac{1} {p} + \frac{{\lambda  +
l\alpha  + l\beta }} {Q} - 1$, then we have%
\begin{equation*}
\left( {\int_G {\left| {Tf\left( u \right)} \right|^q \left| {z_l }
\right|^{ - \alpha q} du} } \right)^{\frac{1} {q}}  \leqslant
C\left( {\int_G {\left| {f\left( u \right)} \right|^p \left| {z'_l }
\right|^{\beta p} du} } \right)^{\frac{1} {p}},
\end{equation*}
where $T$ is the operator in (2.2), $C = C_{\alpha ,\beta
,p,\lambda ,G,l}$ is a positive constant independent of function
$f$, $u = (z_1 ,z_2 , \ldots ,z_r ),v = (z'_1 ,z'_2 , \ldots
,z'_r ) \in G$.

If one takes $f\left( v \right) = \frac{{g\left( v \right)}}
{{\left| {z'_l } \right|^\beta  }}$, then the following equivalent
statement of Theorem 4.1 is evident:

\textbf{Theorem 4.2.} Let $1 < p \leqslant q < \infty ,0 < \lambda
< Q$ and $0 \leqslant \alpha  + \beta  < \frac{{m_l }} {{Q - lm_l
}}\lambda$ such that $\alpha  < m_l /q,\beta  < m_l /p'$ and $
\frac{1} {q} = \frac{1} {p} + \frac{{\lambda  + l\alpha  + l\beta }}
{Q} - 1$, then%
\begin{equation*}
\left\| {\tilde Sg} \right\|_q  \leqslant C\left\| g \right\|_p
\end{equation*}
where $\tilde S$ is an operator of the form%
\begin{equation*}
\tilde Sg\left( u \right) = \tilde S_{\lambda ,\alpha ,\beta }
g\left( u \right) = \int_G {\frac{{g\left( v \right)dv}} {{\left|
{z_l } \right|^\alpha  \left| {u^{ - 1} v} \right|^\lambda  \left|
{z'_l } \right|^\beta  }}},
\end{equation*}
$ C = C_{\alpha ,\beta ,p,\lambda ,G,l}$ is a positive constant
independent of function $f$, $ u = (z_1 ,z_2 , \ldots ,z_r ),v
= (z'_1 ,z'_2 , \ldots ,z'_r ) \in G$

\textbf{Proposition 4.3.} Theorems 2.2 and 4.2 are equivalent.

We easily prove Proposition 4.3 as the proof of Proposition 3.3.
Furthermore, we can give the proof of Theorem 4.1 with the similar
approach in the proof of Theorem 3.1 by using Lemma 3.4 and the
following property.

\textbf{Lemma 4.4.} Suppose $\gamma  < m_l$, then $\left| {z_l }
\right|^{-\gamma}$ is integrable on $\left\{ {\left| u \right| < 1}
\right\}$.

Now we approve the conclusion and list its proof in the next
section.

\textbf{Proof of Theorem 4.1.} Choosing $w_1 \left( u \right) =
\left| {z_l } \right|^{ - \alpha q}$ and $w_2 \left( u \right) =
\left| {z_l } \right|^{\beta p} ,$ it is enough to verify that
(3.5) and (3.6) in Lemma 3.4 are true. As in the proof of
Theorem 3.1, one can easily derive (3.5).

To show (3.6), let us denote%
\begin{eqnarray*}
  M_1 &=&\varphi \left( B \right)\left| B \right|^{\frac{1} {{p'}} + \frac{1}
{q}} , \\
  M_2 &=&\left( {\frac{1} {{\left| B \right|}}\int_B {\left| {z_l } \right|^{
- \alpha q\tau } du} } \right)^{\frac{1} {{q\tau }}} , \\
  M_3 &=&\left( {\frac{1} {{\left| B \right|}}\int_B {\left| {z_l } \right|^{
- \beta p'\tau } du} } \right)^{\frac{1} {{p'\tau }}} ,
\end{eqnarray*}
and estimate $M_i \;(i = 1,2,3)$, respectively.

Setting $\bar \lambda {\text{ = }}Q\left( {\frac{1} {q} + \frac{1}
{{p'}}} \right)$, we have $\frac{1} {q} = \frac{1} {p} + \frac{{\bar
\lambda }} {Q} - 1$ and $\bar \lambda  - \lambda  - l\alpha  -
l\beta  = 0$. A direct computation implies%
\begin{equation*}
M_1  = \varphi \left( B \right)\left| B \right|^{\frac{1} {q} +
\frac{1} {{p'}}}  = Cr^{ - \lambda } r^{Q\left( {\frac{1} {q} +
\frac{1} {{p'}}} \right)}  = Cr^{\bar \lambda  - \lambda };
\end{equation*}
if $\alpha q\tau  < m_l$, we have $\int_{\left| u \right| < 1}
{\left| {z_l } \right|^{ - \alpha q\tau } du}  < C_\tau$ from Lemma
4.4, thus%
\begin{eqnarray*}
  M_2  &=& \left( {\frac{1}
{{\left| B \right|}}\int_{\left| u \right| < r} {\left| {z_l }
\right|^{ - \alpha q\tau } du} } \right)^{\frac{1}
{{q\tau }}}\\
   &=& \left( {\frac{1}
{{\left| B \right|}}\int_{\left| u \right| < 1} {\left| {r^l z_l }
\right|^{ - \alpha q\tau } \left( {d \circ \delta _r } \right)\left(
u \right)} } \right)^{\frac{1}
{{q\tau }}} \\
   &=& \left( {\frac{1}
{{\left| B \right|}}\int_{\left| u \right| < 1} {r^{ - l\alpha q\tau
} \left| {z_l } \right|^{ - \alpha q\tau } \left( {r^Q du} \right)}
} \right)^{\frac{1}
{{q\tau }}}\\
   &=& C_\tau  \left( {\frac{{r^{Q - l\alpha q\tau } }}
{{r^Q }}} \right)^{\frac{1} {{q\tau }}} \left( {\int_{\left| u
\right| < 1} {\left| {z_l } \right|^{ - \alpha q\tau } du} }
\right)^{\frac{1}
{{q\tau }}} \\
   &\leq& C_\tau  r^{ - l\alpha };
\end{eqnarray*}
if $\beta p'\tau  < m_l$, then%
\begin{equation*}
M_3  = \left( {\frac{1} {{\left| B \right|}}\int_B {\left| {z_l }
\right|^{ - \beta p'\tau } du} } \right)^{\frac{1} {{p'\tau }}}
\leqslant C_\tau  r^{ - l\beta }.
\end{equation*}
Since $\alpha  < m_l /q,\beta  < m_l /p',$ we see that there exists $\tau  > 1$
such that $\alpha q\tau  < m_l$ and $\beta p'\tau  < m_l$. Then%
\begin{equation*}
M_1  \cdot M_2  \cdot M_3  \leqslant C_\tau  r^{\bar \lambda  -
\lambda  - l\alpha  - l\beta }  = C_\tau   < \infty,
\end{equation*}
which shows (3.6). Then the proof of Theorem 4.1 is completed by
Lemma 3.4.$\Box$

\section{Proof of Lemma 4.4}

We apply the method in [9, p.383] to check Lemma 4.4.

\textbf{Proof of Lemma 4.4.} Firstly, we consider the case $l =
1\;,\gamma  < m_1$. Thanks to the stratification of $G$ and $ \left|
u \right| = \left( {\sum\limits_{j = 1}^r {\left| {z_j }
\right|^{\frac{{2r!}} {j}} } } \right)^{\frac{1} {{2r!}}}$, it yield%
\begin{eqnarray*}
&&  \int_{\left| u \right| < 1} {\left| {z_1 } \right|} ^{-\gamma}  du \\
   &=& \int_{\left| {z_r } \right| < 1} {dz_r \int_{\left| {z_{r - 1} } \right| < \left( {1 - \left| {z_r } \right|^{\frac{{2r!}}
{r}} } \right)^{\frac{{r - 1}} {{2r!}}} } {dz_{r - 1}   \ldots
\int_{\left| {z_2 } \right| < \left( {1 - \sum\limits_{j = 3}^r
{\left| {z_j } \right|^{^{\frac{{2r!}} {j}} } } } \right)^{\frac{2}
{{2r!}}} } {dz_2 } } }  \cdot \\
  &&\int_{\left| {z_1 } \right| < \left( {1 - \sum\limits_{j = 2}^r {\left| {z_j } \right|^{^{\frac{{2r!}}
{j}} } } } \right)^{\frac{1} {{2r!}}} } {\left| {z_1 }
\right|^{-\gamma}  dz_1 }
\end{eqnarray*}
\begin{eqnarray*}
  &=& \omega _{m_1  - 1} \int_{\left| {z_r } \right| < 1} {dz_r \int_{\left| {z_{r - 1} } \right| < \left( {1 - \left| {z_r } \right|^{\frac{{2r!}}
{r}} } \right)^{\frac{{r - 1}} {{2r!}}} } {dz_{r - 1} \ldots
\int_{\left| {z_2 } \right| < \left( {1 - \sum\limits_{j = 3}^r
{\left| {z_j } \right|^{^{\frac{{2r!}} {j}} } } } \right)^{\frac{2}
{{2r!}}} } {dz_2 } } }  \cdot  \\
  &&\int_{0 < \rho _1  < \left( {1 - \sum\limits_{j = 2}^r {\left| {z_j } \right|^{^{\frac{{2r!}}
{j}} } } } \right)^{\frac{1}
{{2r!}}} } {\rho _1^{-\gamma  + m_1  - 1} d\rho _1 }  \\
   &=& \frac{{\omega _{m_1  - 1} }}
{{-\gamma  + m_1 }}\int_{\left| {z_r } \right| < 1} {dz_r
\int_{\left| {z_{r - 1} } \right| < \left( {1 - \left| {z_r }
\right|^{\frac{{2r!}} {r}} } \right)^{\frac{{r - 1}}
{{2r!}}} } {dz_{r - 1}   \ldots  } } \\
  &&\int_{\left| {z_2 } \right| < \left( {1 - \sum\limits_{j = 3}^r {\left| {z_j } \right|^{^{\frac{{2r!}}
{j}} } } } \right)^{\frac{2} {{2r!}}} } {\left( {1 - \sum\limits_{j
= 2}^r {\left| {z_j } \right|^{^{\frac{{2r!}} {j}} } } }
\right)^{\frac{{-\gamma  + m_1 }} {{2r!}}} dz_2 } ,
\end{eqnarray*}
where $\omega_{n}$ denote the volume of a unit ball in $R^{n}$.

It follows with a computation that%
\begin{eqnarray*}
  &&\int_{\left| {z_2 } \right| < \left( {1 - \sum\limits_{j = 3}^r {\left| {z_j } \right|^{^{\frac{{2r!}}
{j}} } } } \right)^{\frac{2} {{2r!}}} } {\left( {1 - \sum\limits_{j
= 2}^r {\left| {z_j } \right|^{^{\frac{{2r!}} {j}} } } }
\right)^{\frac{{-\gamma  + m_1 }}
{{2r!}}} dz_2 }  \\
   &=& \omega _{m_2  - 1} \int_{0 < \rho _2  < \left( {1 - \sum\limits_{j = 3}^r {\left| {z_j } \right|^{^{\frac{{2r!}}
{j}} } } } \right)^{\frac{2} {{2r!}}} } {\left( {1 - \sum\limits_{j
= 3}^r {\left| {z_j } \right|^{^{\frac{{2r!}} {j}} } }  - \rho
_2^{\frac{{2r!}} {2}} } \right)^{\frac{{-\gamma  + m_1 }}
{{2r!}}} \rho _2^{m_2  - 1} d\rho _2 }  \\
   &=& \omega _{m_2  - 1} \left( {1 - \sum\limits_{j = 3}^r {\left| {z_j } \right|^{^{\frac{{2r!}}
{j}} } } } \right)^{\frac{{-\gamma  + m_1 }}
{{2r!}}}  \cdot  \\
  &&\int_{0 < \rho _2  < \left( {1 - \sum\limits_{j = 3}^r {\left| {z_j } \right|^{^{\frac{{2r!}}
{j}} } } } \right)^{\frac{2} {{2r!}}} } {\left( {1 - \left(
{\frac{{\rho _2 }} {{\left( {1 - \sum\limits_{j = 3}^r {\left| {z_j
} \right|^{^{\frac{{2r!}} {j}} } } } \right)^{\frac{2} {{2r!}}} }}}
\right)^{\frac{{2r!}} {2}} } \right)^{\frac{{-\gamma  + m_1 }}
{{2r!}}} \rho _2^{m_2  - 1} d\rho _2 }  \\
   &=& \omega _{m_2  - 1} \left( {1 - \sum\limits_{j = 3}^r {\left| {z_j } \right|^{^{\frac{{2r!}}
{j}} } } } \right)^{\frac{{-\gamma  + m_1  + 2m_2  - 2}} {{2r!}}}
\int_{0 < \rho _2  < 1} {\left( {1 - \rho _2^{\frac{{2r!}} {2}} }
\right)^{\frac{{-\gamma  + m_1 }} {{2r!}}} \rho _2^{m_2  - 1} d\rho
_2 }\\
  & =& \omega _{m_2  - 1} \frac{2}
{{2r!}}\left( {1 - \sum\limits_{j = 3}^r {\left| {z_j }
\right|^{^{\frac{{2r!}} {j}} } } } \right)^{\frac{{-\gamma  + m_1  +
2m_2  - 2}} {{2r!}}} \int_{0 < \rho _2  < 1} {\left( {1 - \rho _2 }
\right)^{\frac{{-\gamma  + m_1 }} {{2r!}}} \rho _2^{\frac{{2m_2 }}
{{2r!}} - 1} d\rho _2 }  \\
   &=& \omega _{m_2  - 1} \frac{2}
{{2r!}}{\rm B}\left( {\frac{{-\gamma  + m_1 }} {{2r!}} +
1,\frac{{2m_2 }} {{2r!}}} \right)\left( {1 - \sum\limits_{j = 3}^r
{\left| {z_j } \right|^{^{\frac{{2r!}} {j}} } } }
\right)^{\frac{{-\gamma  + m_1  + 2m_2  - 2}} {{2r!}}} ,
\end{eqnarray*}
where ${\rm B}\left( { \cdot , \cdot } \right)$ denotes a
Beta-function. By the induction, we obtain%
\begin{eqnarray*}
 && \int_{\left| u \right| < 1} {\left| {z_1 } \right|} ^{-\gamma}  du \hfill \\
   &=& \left( {\prod\limits_{j = 1}^r {\omega _{m_j  - 1} } } \right)\frac{{r!}}
{{\left( {2r!} \right)^{r - 1} }}\frac{1} {{-\gamma  + m_1 }}\left(
{\prod\limits_{j = 2}^r {{\rm B}\left( {\frac{{-\gamma  + 1 +
\sum\limits_{i = 1}^{j - 1} {i\left( {m_i  - 1} \right)} }} {{2r!}}
+ 1,\frac{{jm_j }} {{2r!}}} \right)} } \right).
\end{eqnarray*}

Secondly, if $\gamma  < m_l \left( {l = 2, \ldots ,r} \right)$, then
we have%
\begin{eqnarray*}
  &&\int_{\left| u \right| < 1} {\left| {z_l } \right|} ^{-\gamma}  du \\
   &=& \int_{\left| {z_r } \right| < 1} {dz_r   \ldots   \int_{\left| {z_{l + 1} } \right| < \left( {1 - \sum\limits_{j = l + 2}^r {\left| {z_j } \right|^{^{\frac{{2r!}}
{j}} } } } \right)^{\frac{{l + 1}}
{{2r!}}} } {dz_{l + 1}  \cdot } }  \\
  &&\int_{\left| {z_{l - 1} } \right| < \left( {1 - \sum\limits_{j = l + 1}^r {\left| {z_j } \right|^{^{\frac{{2r!}}
{j}} } } } \right)^{\frac{{l - 1}} {{2r!}}} } {dz_{l - 1} }
\int_{\left| {z_{l - 2} } \right| < \left( {1 - \sum\limits_{j = l -
1,j \ne l}^r {\left| {z_j } \right|^{^{\frac{{2r!}} {j}} } } }
\right)^{\frac{{l - 2}}
{{2r!}}} } {dz_{l - 2}    \ldots }    \\
  &&\int_{\left| {z_1 } \right| < \left( {1 - \sum\limits_{j = 2,j \ne l}^r {\left| {z_j } \right|^{^{\frac{{2r!}}
{j}} } } } \right)^{\frac{1} {{2r!}}} } {dz_1 \int_{\left| {z_l }
\right| < \left( {1 - \sum\limits_{j = 1,j \ne l}^r {\left| {z_j }
\right|^{^{\frac{{2r!}} {j}} } } } \right)^{\frac{l}
{{2r!}}} } {\left| {z_l } \right|^{-\gamma}  } dz_l }  \\
   &=& \omega _{m_l  - 1} \int_{\left| {z_r } \right| < 1} {dz_r    \ldots   \int_{\left| {z_{l + 1} } \right| < \left( {1 - \sum\limits_{j = l + 2}^r {\left| {z_j } \right|^{^{\frac{{2r!}}
{j}} } } } \right)^{\frac{{l + 1}}
{{2r!}}} } {dz_{l + 1}  \cdot } }  \\
  &&\int_{\left| {z_{l - 1} } \right| < \left( {1 - \sum\limits_{j = l + 1}^r {\left| {z_j } \right|^{^{\frac{{2r!}}
{j}} } } } \right)^{\frac{{l - 1}} {{2r!}}} } {dz_{l - 1} }
\int_{\left| {z_{l - 2} } \right| < \left( {1 - \sum\limits_{j = l -
1,j \ne l}^r {\left| {z_j } \right|^{^{\frac{{2r!}} {j}} } } }
\right)^{\frac{{l - 2}}
{{2r!}}} } {dz_{l - 2}    \ldots }    \\
  &&\int_{\left| {z_1 } \right| < \left( {1 - \sum\limits_{j = 2,j \ne l}^r {\left| {z_j } \right|^{^{\frac{{2r!}}
{j}} } } } \right)^{\frac{1} {{2r!}}} } {dz_1 \int_{0 < \rho _l  <
\left( {1 - \sum\limits_{j = 1,j \ne l}^r {\left| {z_j }
\right|^{^{\frac{{2r!}} {j}} } } } \right)^{\frac{l}
{{2r!}}} } {\rho _l^{-\gamma  + m_l  - 1} d\rho _l } }  \\
   &=& \frac{{\omega _{m_l  - 1} }}
{{-\gamma  + m_l }}\int_{\left| {z_r } \right| < 1} {dz_r \ldots
 \int_{\left| {z_{l + 1} } \right| < \left( {1 -
\sum\limits_{j = l + 2}^r {\left| {z_j } \right|^{^{\frac{{2r!}}
{j}} } } } \right)^{\frac{{l + 1}}
{{2r!}}} } {dz_{l + 1}  \cdot } }  \\
  &&\int_{\left| {z_{l - 1} } \right| < \left( {1 - \sum\limits_{j = l + 1}^r {\left| {z_j } \right|^{^{\frac{{2r!}}
{j}} } } } \right)^{\frac{{l - 1}} {{2r!}}} } {dz_{l - 1}   }
\int_{\left| {z_{l - 2} } \right| < \left( {1 - \sum\limits_{j = l -
1,j \ne l}^r {\left| {z_j } \right|^{^{\frac{{2r!}} {j}} } } }
\right)^{\frac{{l - 2}}
{{2r!}}} } {dz_{l - 2}   \ldots }    \\
  &&\int_{\left| {z_1 } \right| < \left( {1 - \sum\limits_{j = 2,j \ne l}^r {\left| {z_j } \right|^{^{\frac{{2r!}}
{j}} } } } \right)^{\frac{1} {{2r!}}} } {\left( {1 - \sum\limits_{j
= 1,j \ne l}^r {\left| {z_j } \right|^{^{\frac{{2r!}} {j}} } } }
\right)^{\frac{{-l\gamma  + lm_l }} {{2r!}}} dz_1 }
\end{eqnarray*}

Similarly to the proof of $l = 1$, it implies%
\begin{eqnarray*}
   &&\int_{\left| u \right| < 1} {\left| {z_l } \right|} ^{-\gamma}  du\\
  &=& \left( {\prod\limits_{j = 1}^r {\omega _{m_j  - 1} } } \right)
\frac{{r!}} {{l\left( {2r!} \right)^{r - 1} }}\frac{1} {{-\gamma  +
m_l }}\left( {\prod\limits_{j = 1,j \ne l}^r {{\rm B}\left( {\frac{{
- l\gamma  + lm_l  + \sum\limits_{i = 1,i \ne l}^{j - 1} {i\left(
{m_i  - 1} \right)} }} {{2r!}} + 1,\frac{{jm_j }} {{2r!}}} \right)}
} \right).
\end{eqnarray*}
The proof is completed.$\Box$

\textbf{Acknowledgment}
This work was supported by the National Natural
Science Foundation of China (Grant No. 11271299), the Mathematical
Tianyuan Foundation of China(Grant No. 11126027) and Natural Science
Foundation Research Project of Shaanxi Province
(2012JM1014); *corresponding author.


\end{document}